\documentclass[12pt]{article}
\usepackage[utf8x]{inputenc}
\usepackage{amssymb}
\usepackage{color}
\usepackage{geometry}

\def\1{\mathbf{1}}

\newif\ifhide
\hidetrue

\newtheorem{lem} {Lemma}
\newtheorem{thm} {Theorem}

\newtheorem{rem} {Remark}
\newtheorem{cor} {Corollary}

\newif\ifmoment
\momentfalse

\title{{\normalsize\tt\hfill\jobname.tex}\\
\bf On recurrence and availability factor  
for single--server system with general arrivals
}
\author{A. Yu. Veretennikov\footnote{ University of Leeds, UK; National Research University Higher School of Economics, and Institute for Information Transmission Problems, Moscow, Russia, email: a.veretennikov @ leeds.ac.uk. The work was prepared within the framework of a subsidy granted to the HSE by the Government of the Russian Federation for the implementation of the Global Competitiveness Program, and supported by the RFBR grant 14-01-00319-a.}
}
\begin{document}
\maketitle

\begin{abstract}
Recurrence and ergodic properties are established for a single--server queueing system with variable intensities of arrivals and service. Convergence to stationarity is also interpreted in terms of reliability theory.
\end{abstract}

\section{Introduction}
In the last decades, queueing systems generalising  $M/G/1/\infty$, or
$M/G/1$ (cf. \cite{GK}) -- one of the  most  important queueing systems -- 
attracted much attention, see \cite{Asmussen} -- \cite{fakinos1987},
%\cite{Asmussen}, \cite{Bambos89},  \cite{borovkov2003}, \cite{bramson2008}, \cite{fakinos1987},
 \cite{Thor83}.
%\cite{Kalyanaraman2008}, \cite{masuyama2003},
%\cite{Miyoshi2001},
%\cite{osipova2009}, \cite{schonlein2012}?]. 
In this paper a single--server system similar to \cite{Ve2013_ait, Ver_Z-2014}
is considered,
%{\color{red}(cf. also \cite{Thor00})},
in which {\em intensities} of new arrivals as well as of their
service may depend on the ``whole state'' of the system and the
whole state includes the number of customers in the system --
waiting and on service -- {\em and} on the elapsed time of the last
service, as well as on the elapsed time since the end of the last service. Batch arrivals are not allowed. The news in comparison to \cite{Ve2013_ait, Ver_Z-2014} is that at any state, even if the system idle (no service), the intensity of new arrivals  
may depend on the time from the last end of service. The details of the system description will be formalised in the beginning of the next section.
By the {\em m-availability factor} of the system we understand the probability of the idle state 
if $m=0$, or probability of $m$ customers in total on the server and in the queue. We do not use notation $G/G/1$ (or $GI/GI/1$) only because some conditions on intensities are assumed, which makes the model  slightly less general. 
The problem addressed in the paper is how to estimate convergence rate of characteristics of the system including the $m$-availability factors 
to their stationary values. 

~

The {\em
elapsed} service time is assumed to be known at any moment, but 
the remaining service times for each customer are not. For definiteness,
the discipline of serving is FIFO, although other disciplines may be also considered.  

~

The paper consists of the Section 1 -- Introduction, of the setting and main result in the Section 2, of the auxiliary lemmata in the Section 3 and of the short sketch of the proof of the main result in the Section 4.

\section{The setting and main results}

\subsection{Defining the process}

Let us present the class of models under investigation in this
paper. Here the state space is a union of subspaces,
\[
{\mathcal X} = 
%\{(0,x), \, x\ge 0)\} 
\{(0,y):  \;
y\ge 0\}
\cup
\bigcup_{n=1}^{\infty} \{(n,x,y):  \;
x,y\ge 0\}.
\]
% with topology arising from the metric
%($X=(n,x)$, $X'=(n',x')$)
%\[
%\mbox{dist}(X, X')=|x-x'|\,1(n'=n) + |n-n'|.
%\]
Functions of
class $C^1({\mathcal X})$ are understood as functions with
classical continuous derivatives with respect to the variable $x$.
Functions with compact support  on ${\mathcal X}$ are understood
as functions vanishing outside some domain bounded in this metric:
for example,  $C^1_0({\mathcal X})$ stands for the class of
functions with compact support and one continuous derivative.
There is a generalised Poisson arrival flow with intensity $
\lambda(X), $ where \( X = (n,x,y) \; \mbox{for any $n \ge 1$} \), and \( X = (0,y) \; \mbox{for $n = 0 $}\). Slightly abusing notations, it is convenient to write $X=(n,x,y)$ for $n=0$ as well, assuming that in this case $x=0$. If $n>0$, then the
server is serving one customer while all others are waiting in a
queue. When the last service ends, immediately a new service of  the next customer from the queue starts. If $n=0$ then the server remains idle until the next customer 
arrival; the intensity of such arrival at state $(0,y)\equiv (0,0,y)$ may be variable depending on the value $y$, which stands for the elapsed time from the last end of service. Here $n$ denotes the total number of customers in the
system, and $x$ stands for the elapsed time of the current service  (except for $n=0$, which was explained earlier), and $y$ is the elapsed from from the last arrival. {\em Normally}, intensity of arrivals depend on $n$ and $y$, while intensity of service depends on $n$ and $x$; however, we allow more general dependence. 
Denote $n_t=n(X_t)$ -- the number of customers corresponding to
the state $X_t$, and $x_t=x(X_t)$, the second component of the
process $(X_t)$, and $y_t=y(X_t)$, the third component of the
process $(X_t)$ (the third if $n>0$)).  For
any $X=(n,x,y)$, intensity of service $h(X) \equiv h(n,x,y)$ is
defined; it is also convenient to assume that $h(X)=0$ for $n(X)=0$.
Both intensities $\lambda$ and $h$ are understood in the following
way, which is a definition: on any nonrandom interval of
time $[t,t+\Delta)$, conditional probability given $X_t$ that the
current service will {\em not} be finished and there will be no
new arrivals reads,
\begin{equation}\label{intu1}
\exp\left(-\int_0^{\Delta} (\lambda+h)
(n_{t},x_{t}+s, y_t+s)\,ds\right).
\end{equation}
In the sequel, $\lambda$ and $h$ are assumed to be {\em bounded}.
In this case, for $\Delta>0$ small enough, the expression in
(\ref{intu1}) may be rewritten as
\begin{equation}\label{intu2}
1-\int_0^{\Delta} (\lambda+h)(n_{t},x_{t}+s, y_{t}+s)\,ds + O(\Delta^2),
\qquad \Delta\to 0,
\end{equation}
and this what is ``usually'' replaced by
\[
1 -  (\lambda(X_t)+h(X_t))\Delta + O(\Delta^2).
\]
However, in our situation, the latter replacement may be incorrect because
of discontinuities of the functions $\lambda$ and $h$. Emphasize
that from time $t$ and until the next jump, the evolution of the
process $X$ is {\em deterministic}, which makes the process {\em piecewise-linear Markov,} see, e.g., \cite{GK}. The (conditional given $X_t$)
density  of the moment of a new arrival {\em or} of the end of the
current service after $t$ at $x_t+ z$, $z\ge 0$ equals,
\begin{equation}\label{intu33}
 (\lambda(n_t, x_t+z, y_t+z)+h(n_t, x_t+z, y_t+z))
 \exp\left(-\int_0^{\Delta}
(\lambda+h)(n_t, x_{t}+s, y_t+s))\,ds\right).
\end{equation}
Further, given $X_t$, the moments of the next ``candidates'' for
jumps up and down are conditionally independent and have the
(conditional -- given $X_t$) density, respectively,
\begin{equation}\label{intu5}
\begin{array}{c}
   \lambda(X_t+z)\exp\left(-\int\limits_0^{z}
\lambda(X_{t}+s)\,ds\right) \; \\
   \\
  \mbox{and} \; \\
   \\
  h(X_t+z)\exp\left(-\int\limits_0^{z} h(X_{t}+s)\,ds\right), \; z\ge 0.\\
\end{array}
\end{equation}
(Here $X_{t}+s := (n_t, x_t+s, y_t+s)$.) Notice that (\ref{intu33}) does correspond to conditionally
independent densities given in (\ref{intu5}).

\subsection{Main result}
Let 
\[
\Lambda:= \sup_{n,x,y: \,n>0}\lambda(n,x,y) < \infty.
\]
For establishing convergence rate to the stationary regime, we
assume similarly to \cite{Ve2013_ait, Ver_Z-2014},
\begin{equation}\label{eq2}
\inf_{n>0, y} h(n,x,y) \ge \frac{C_0}{1+x}, \quad x\ge 0. 
\end{equation}
We also assume a new condition related to $\lambda_0(t) = \lambda(0, 0,t)$, which was   constant in the earlier papers:  
now it is allowed to be variable and satisfying 
\begin{equation}\label{eq20}
0<\inf_{t\ge 0}\lambda_0(t) \le \sup_{t\ge 0}\lambda_0(t) <\infty. 
\end{equation}
Recall that the process has no explosion with probability one due
to the boundedness of both intensities, i.e., the trajectory may
have only finitely many jumps on any finite interval of time.

\begin{thm}\label{thm2}
Let the functions $\lambda$ and $h$ be Borel measurable and
bounded and let the assumptions (\ref{eq2}) and (\ref{eq20}) be satisfied. Then, under the assumptions above, if $C_0$ is large enough, then 
there exists a unique stationary measure $\mu$. Moreover, for any $m>k$, $C>0$ there exists $\bar C>0$ such that if $C_0\ge \bar C$, then for any $t\ge 0$,
\begin{equation}\label{est}
 \|\mu^{n,x,y}_t - \mu \|_{TV} \le
 C \,\frac{(1+n+x+y)^m}{(1+t)^{k+1}},
\end{equation}
where $\mu^{n,x,y}_t$ is a marginal distribution of the process
$(X_t, \, t\ge 0)$ with the initial data $X=(n,x,y)\in {\mathcal
X}$.
\end{thm}

\begin{rem}
It is plausible that the bound in (\ref{est}) may be improved so that the right hand side does not depend on $y$. 
Moreover, given all other constants, the value $C$ in (\ref{est}) may be made ``computable'', with a rather involved but explicit dependence on other constants. Moreover, it is likely that the condition (\ref{eq20}) may be replaced by a weaker one, 
\begin{equation}\label{eq21}
\frac{C_0'}{1+t} \le \lambda_0(t) \le \sup_{t\ge 0}\lambda_0(t) <\infty, 
\end{equation}
along with the assumption that $C'_0$ is large enough. 
However, all  these issues require a bit more accuracy in the calculus and we do not pursue these goals here leaving them until further publications with complete technical details.
\end{rem}

\section{Lemmata}
Recall \cite{Dynkin} that the generator of a Markov process $(X_t,
\, t\ge 0)$ is an operator \( {\mathcal G}, \) such that for a
sufficiently large class of functions $f$
\begin{equation}\label{dynkin2}
 \sup_X \lim_{t\to 0} \left\|\frac{E_Xf(X_t) - f(X)}{t} -
{\mathcal G}f(X)\right\| = 0
\end{equation}
in the norm of the state space of the process; the notion of
generator does depend on this norm.
An
operator ${\mathcal G}$ is called a {\em mild generalised
generator} (another name is extended generator) if (\ref{dynkin2}) is replaced by its corollary
(\ref{dynkin1}) below called {\em Dynkin's formula},
or {\em Dynkin's identity}  \cite[Ch. 1, \S 3]{Dynkin},
\begin{equation}\label{dynkin1}
E_Xf(X_t) - f(X) = E_X\int_0^t {\mathcal G}f(X_s)\,ds,
\end{equation}
also for a wide enough class of functions $f$. 
We will also use the 
non-homogeneous counterpart of Dynkin's formula,  
\begin{equation}\label{dynkin_t}
E_X\varphi(t,X_t) - \varphi(0,X) = E_X\int_0^t
\left(\frac{\partial}{\partial s}\varphi(s, X_s) + {\mathcal
G}\varphi(s,X_s)\right)\,ds,
\end{equation}
for appropriate functions of two variables $(\varphi(t,X))$. Both (\ref{dynkin1}) and (\ref{dynkin_t}) play
a very important role in analysis of Markov models and under our assumptions may be justified similarly to \cite{Ver_Z-2014}. Here
 $X$ is a (non-random)
initial value of the process. 
Both formulae (\ref{dynkin1})--(\ref{dynkin_t}) hold true for a
large class of functions $f$, $\varphi$ with ${\mathcal G}$ 
given by the standard expression,

\begin{eqnarray*}
{\mathcal G}f(X) := \frac{\partial}{\partial x}f(X)1(n(X)>0) + \frac{\partial}{\partial y}f(X) 
 \\\\
+  \lambda(X)(f(X^+) - f(X))
+ h(X) (f(X^-) - f(X)),
\end{eqnarray*}
where for any $X=(n,x,y)$,
\[
X^+ := (n+1,x,0), \quad X^-:= ((n-1)\vee 0,0,y)
\]
(here $a\vee b = \max(a,b)$).
Under our 
minimal assumptions on regularity of intensities this may be justified similarly to \cite{Ver_Z-2014}.

\begin{lem}\label{thm1}
If the functions $\lambda$ and $h$ are Borel measurable and
bounded, then the formulae (\ref{dynkin1}) and (\ref{dynkin_t})
hold true for any $t>0$ for every $f\in C^1_b({\mathcal X})$ and
$\varphi\in C^1_b([0,\infty)\times {\mathcal X})$, respectively.
%\end{thm}
Moreover, the process $(X_t, \, t\ge 0)$ is strong Markov with
respect to the filtration \(({\mathcal F}^X_t, \, t\ge 0)\).
\end{lem}

~

\noindent
Further, let
\begin{equation}\label{LL}
L_m(X) = (n+1+x +y)^m, \quad L_{k,m}(t,X) = (1+t)^k L_{m}(X).
\end{equation}
The extensions of Dynkin's formulae for some  unbounded functions hold true: we will need them for the Lyapunov functions in (\ref{LL}). 
\begin{cor}\label{cor1}
Under the assumptions of the Lemma \ref{thm1},
\begin{eqnarray}\label{M2}
L_{m}(X_t) - L_{m}(X) = \int_0^t \lambda(X_s)
\left[\phantom{\!\!\!\frac{}{} } \left(L_{m}(X^{(+)}_s) -
L_{m}(X_s)\right) \right.
\hspace{2cm}
 \nonumber \\ \\ \nonumber
\left. + h(X_s) \left(L_{m}(X^{-}_s) - L_{m} (X_s)\right) 
+1(n(X_s)>0) \frac{\partial}{\partial x}L_{m}(X_s)
+\frac{\partial}{\partial y}L_{m}(X_s)\right]\,ds +M_t,
 \end{eqnarray}
with some martingale $M_t$, and also
\begin{eqnarray}\label{M2t}
L_{k,m}(t,X_t) - L_{k,m}(0,X) = \int_0^t \left[\lambda(X_s)
\left(L_{k,m}(s,X^{(+)}_s) - L_{k,m}(s,X_s)\right) \right.
\hspace{2cm}
  \nonumber \\ \\ \nonumber \left.
+ h(X_s) \left(L_{k,m}(s,X^{-}_s) - L_{k,m}(s,X_s)\right) +
\left(1(n(X_s)>0)\frac{\partial}{\partial x} + \frac{\partial}{\partial y}+ \frac{\partial}{\partial
s}\right) L_{k,m}(s,X_s)\right]\,ds +\tilde M_t,
 \end{eqnarray}
with some martingale
$\tilde M_t$.
%\end{lem}
\end{cor}
About a martingale approach in queueing models  see, for example, \cite{LSh}. The proof of the Lemma \ref{thm1} is based on the next three Lemmata. 
The first of them is a rigorous statement concerning a 
well-known folklore property that probability of ``one event'' on a small
nonrandom interval of length $\Delta$ is of the order $O(\Delta)$
and probability of ``two or more  events'' on the same interval is
of the order $O(\Delta^2)$. Of course, this is a common knowledge in queueing theory, yet  for
discontinuous intensities it has to be, at least, explicitly stated.
% justified. 

\begin{lem}\label{lem1}
Under the assumptions of the Theorem \ref{thm1}, for any $t\ge 0$,
\begin{equation}\label{z0}
P_{X_{t}}(\mbox{no jumps on $(t, t+\Delta]$}) =
\exp(-\int_0^\Delta (\lambda+h)(X_t+s)\,ds) \quad (= 1 +
O(\Delta)),
\end{equation}
\begin{equation}\label{z1}
P_{X_{t}}(\mbox{at least one jump on $(t, t+\Delta]$}) =
O(\Delta),
\end{equation}
\begin{equation}\label{z1up}
P_{X_{t}}(\mbox{exactly one jump up \&  no down on $(t,
t+\Delta]$}) = \int_0^\Delta \lambda(X_t+s)\,ds + O(\Delta^2),
\end{equation}
\begin{equation}\label{z1down}
P_{X_{t}}(\mbox{exactly one jump down \& no up on $(t,
t+\Delta]$}) = \int_0^\Delta h(X_t+s)\,ds + O(\Delta^2),
\end{equation}
and
\begin{equation}\label{z2}
P_{X_{t}}(\mbox{at least two jumps on $(t, t+\Delta]$}) =
O(\Delta^2).
\end{equation}
In all cases above, $O(\Delta)$ and  $O(\Delta^2)$ are uniform
with respect to $X_{t}$ and only depend on the norm
$\sup_{X}(\lambda(X)+h(X))$, that is, there exist $C>0, \,
\Delta_0>0$ such that for any $X$ and any $\Delta<\Delta_0$,
\begin{eqnarray}\label{z_uni}
\limsup_{\Delta\to 0} \left\{\Delta^{-1}P_{X}(\mbox{at least one
jumps on $(0,\Delta]$})\right.
 %\\\\
+ \Delta^{-2} P_{X}(\mbox{at least two jumps on $(0, \Delta]$})
 \nonumber \\ \nonumber \\ \nonumber
+ \Delta^{-2}\left[P_{X_{t}}(\mbox{one jump up \& no down on $(t,
t+\Delta]$}) - \int_0^\Delta \lambda(X_t+s)\,ds \right]
 \nonumber \\  \\ \nonumber
\left.+ \Delta^{-2}\left[P_{X_{t}}(\mbox{one jump down \& no up on
$(t, t+\Delta]$}) - \int_0^\Delta h(X_t+s)\,ds \right]\right\}
<C<\infty.
% \nonumber \\ \nonumber \\
%\le C. \hspace{6cm}
\end{eqnarray}
\end{lem}
The next two Lemmata are needed for the justification that the process with discontinuous intensities is, indeed, strong Markov.

\begin{lem}\label{Le2}
Under the assumptions of the Theorem \ref{thm1}, the semigroup
\linebreak $T_tf(X) = E_Xf(X_t)$ is continuous in $t$.
\end{lem}

%\noindent Denote $X'=(n',x')\uparrow X=(n,x)$ iff $x'\uparrow x$
%and $n' = n$ when $x'$ is close enough to $x$. Similarly,
%$X'=(n',x')\downarrow X=(n,x)$ iff $x'\downarrow x$ and $n' = n$
%for $x'$ close enough to $x$.

\begin{lem}\label{lem4}
Under the assumptions of the Theorem \ref{thm1} the process
$(X_t,\,t\ge 0)$ is Feller, that is, $T_tf(\cdot)\in C_b({\mathcal
X})$ for any $f\in C_b({\mathcal X})$.
\end{lem}

\medskip

\noindent
The proofs of all Lemmata may be performed similarly to \cite{Ver_Z-2014}.

\section{Sketch of Proof of Theorem \ref{thm2}}
The proof of convergence in total variation with rate of convergence repeats the calculus in \cite{Ve2013_ait} based on the Lyapunov functions \(L_{m}(X)\)  and \(L_{k,m}(t,X)\) from (\ref{LL}), 
and on Dynkin's formulae (\ref{dynkin1}) and (\ref{dynkin_t}) due
to the Corollary \ref{cor1}. Without big changes, this calculus provides  a polynomial moment bound 
\begin{equation}\label{pb}
E_X \tau_0^k \le C L_{m}(X) \le C (n+1 + x + y)^m,
\end{equation}
for certain values of $k$ and for the hitting time
\[
\tau_0:= \inf(t\ge 0: \; n_t = 0).
\]
Namely, once the process attains the set  $\{n=0\}$, it may be successfully coupled with another (stationary) version of the same process at their joint jump $\{n=0\} \, \mapsto \{n=1\}$. This is because, in particular, immediately after such a jump the state of each process reads as $(1, 0, 0)$; in other words, this is a regeneration state.
The news is only a wider class of
intensities, which may be all variable (as well as  discontinuous) including $\lambda_0$; however, this affects the calculus only a little, once it is established that (\ref{dynkin1}) and
(\ref{dynkin_t}) hold true, because this calculus involves only time values $t<\tau_0$. (Some change will be in the procedure of coupling, though.) In turn, the inequality (\ref{pb}) provides a bound for the rate of convergence, for the justification of which rate there are various approaches such as versions of coupling as well as  renewal theory. Convergence of probabilities in the definition of $m$-availability factors is a special case of a more general convergence  in total variation. We drop further details, which will be specified in a further publication.

\end{document}